\documentclass[12pt,reqno]{amsart}
\usepackage{color}
\usepackage[english]{babel}
\usepackage{amsfonts}
\usepackage{amsmath}
\usepackage{amsthm}
\usepackage{amssymb}
\usepackage[misc]{ifsym}
\usepackage{bbm}
\usepackage{mathrsfs}
\usepackage{esint}
\usepackage{faktor}
\usepackage[utf8]{inputenc}
\usepackage{afterpage}
\usepackage[left=2.9cm,right=2.9cm,top=2.8cm,bottom=2.8cm]{geometry}
\usepackage{graphicx}
\usepackage{enumitem}
\usepackage[dvipsnames]{xcolor}
\usepackage[colorlinks=true,urlcolor=blue,citecolor=blue,linkcolor=blue,linktocpage,pdfpagelabels,bookmarksnumbered,bookmarksopen]{hyperref}

\newcommand{\N}{\mathbb{N}}

\newcommand{\R}{\mathbb{R}}

\newcommand{\D}{\mathcal{D}}
\newcommand{\Div}{\mathrm{div} \, }
\newcommand{\dx}{\, {\rm d} x}

\newcommand{\ds}{\, {\rm d} s}

\newcommand{\drho}{\, {\rm d} \rho}
\newcommand{\eps}{\varepsilon}
\newcommand{\loc}{{\rm loc}}

\DeclareMathOperator{\dist}{dist}

\newtheorem{lemma}{Lemma}[section]
\newtheorem{thm}[lemma]{Theorem}

\theoremstyle{definition}

\newtheorem{rmk}[lemma]{Remark}

\numberwithin{equation}{section}

\begin{document}

\title[Decay estimates for solutions to non-autonomous critical problems]{Decay estimates for solutions to \\ non-autonomous critical $p$-Laplace problems}

\author[L. Baldelli]{Laura Baldelli}
\address[L. Baldelli]{IMAG, Departamento de Análisis Matemático, Universidad de Granada, Campus Fuentenueva, 18071 Granada, Spain}
\email{labaldelli@ugr.es}

\author[U. Guarnotta]{Umberto Guarnotta}
\address[U. Guarnotta]{Dipartimento di Ingegneria Industriale e Scienze Matematiche, Università Politecnica delle Marche, Via Brecce Bianche 12, 60131 Ancona, Italy}
\email{u.guarnotta@univpm.it}

\begin{abstract}
We prove optimal decay estimates for positive solutions to elliptic $p$-Laplacian problems in the entire Euclidean space, when a critical nonlinearity with a decaying source term is considered. Also gradient decay estimates are furnished. Our results extend previous theorems in the literature, in which a purely critical reaction is treated. The technique is based on a priori estimates, regularity results, and rescaling arguments, combined with the doubling lemma.
\end{abstract}

\maketitle

{
\let\thefootnote\relax
\footnote{\Letter \quad Corresponding author: Laura Baldelli.}
}
\setcounter{footnote}{0}

\section{Introduction}

We study the decay of any $u\in \D^{1,p}_0(\R^N)$ solution to
\begin{equation}\label{prob}\tag{P}\begin{cases}
0\le -\Delta_p u \le f(x) + \Lambda u^{p^*-1} &\quad \mbox{in}\;\; \R^N,\\
 u>0 &\quad \mbox{in}\;\; \R^N,
\end{cases}\end{equation}
where $1<p<N$, $\Delta_p u:=\Div(|\nabla u|^{p-2}\nabla u)$ is the $p$-Laplacian operator, and $f:\R^N\to\R$ is a measurable function satisfying
\begin{enumerate}[label={$({\rm H})$},ref={H}]
\item \label{hyp} there exist $\alpha>N$ and $L>0$ such that
$$ 0\le f(x)\le L(1+|x|)^{-\alpha} \quad \mbox{for almost all} \;\; x\in\R^N. $$
\end{enumerate}
Here, with $\D^{1,p}_0(\R^N)$ we denote the $p$-th Beppo Levi space (see Section 2 for details). The differential inequalities in \eqref{prob} are understood in weak sense, namely,
$$ 0\leq \int_{\R^N} |\nabla u|^{p-2}\nabla u\nabla \phi \dx \leq \int_{\R^N} f\phi \dx + \Lambda \int_{\R^N} u^{p^*-1}\phi \dx \quad \mbox{for all} \;\; \phi\in \D^{1,p}_0(\R^N)_+. $$

We provide double-sided decay estimates for both the solutions to \eqref{prob} and their gradients, as stated in the main result below. Apropos, the optimality of these estimates (in terms of decay rate at infinity) is attained by Talenti's functions
$$u_{\varepsilon,x_0}(x)=c_{p,N}\left[\frac{\varepsilon^{\frac{1}{p-1}}}{\varepsilon^{\frac{p}{p-1}}+|x-x_0|^{\frac{p}{p-1}}}\right]^{\frac{N-p}{p}},\qquad c_{p,N}=\left(N^{\frac{1}{p}}\left(\frac{N-p}{p-1}\right)^{\frac{p-1}{p}}\right)^{\frac{N-p}{p}},$$
where $\varepsilon>0$ and $x_0\in\R^N$, which solve the critical equation $-\Delta_p u = u^{p^*-1}$ in $\R^N$.

\begin{thm}\label{mainthm}
Let $1<p<N$, $f:\R^N\to\R$ be a measurable function fulfilling \eqref{hyp}, and $u\in\D^{1,p}_0(\R^N)$ be a solution to \eqref{prob}. Then $u\in C^{1,\tau}_\loc(\R^N)$ for some $\tau\in(0,1)$ and there exist $C_0>0$, depending on $p,N,u$, and $C_1,C_2,C_3,R>0$, depending on $N,p,\Lambda,\alpha,u$, such that 
\begin{equation}\label{main1}
C_0\left(1+|x|^{\frac{N-p}{p-1}}\right)^{-1} \leq u(x)\leq C_1\left(1+|x|^{\frac{N-p}{p-1}}\right)^{-1}
\end{equation}
and
\begin{equation}\label{main2}
\left|\nabla u(x)\right|\leq C_2\left(1+\left|x\right|^{\frac{N-1}{p-1}}\right)^{-1}
\end{equation}
for all $x\in\R^N$, as well as
\begin{equation}\label{main3}
\left|\nabla u(x)\right| \geq C_3|x|^{\frac{1-N}{p-1}} \quad \mbox{for all} \;\; x\in B_R^e. 
\end{equation}
\end{thm}

\begin{rmk}
\label{refteo}
The general nature of the problem \eqref{prob} allows to encompass a broader range of nonlinearities, i.e.,
$$ 0 \leq -\Delta_p u \leq \hat{f}(x) + V(x)u^{k-1} + \hat{\Lambda} u^{p^*-1} \quad \mbox{in} \;\; \R^N, $$
where $\hat{\Lambda}>0$, $\hat{f}$ fulfills \eqref{hyp}, $k\in (1,p^*)$, and $V$ is a non-negative measurable function satisfying, for some $\beta>N\,\frac{p^*-k}{p^*-1}$ and $\hat{L}>0$,
$$ V(x) \leq \hat{L}(1+|x|)^{-\beta} \quad \mbox{for almost all} \;\; x\in\R^N. $$
Indeed, Young's inequality entails
$$ Vu^{k-1} \leq \frac{p^*-k}{p^*-1}\,V^{\frac{p^*-1}{p^*-k}}+\frac{k-1}{p^*-1}\,u^{p^*-1} \quad \mbox{in} \;\; \R^N, $$
so one can choose $f:=\hat{f}+\frac{p^*-k}{p^*-1}\,V^{\frac{p^*-1}{p^*-k}}$ and $\Lambda:=\hat{\Lambda}+\frac{k-1}{p^*-1}$ in \eqref{prob}. Notice that the limit cases $k=1$ and $k=p^*$ correspond, respectively, to $\hat{f}(x)$ and $u^{p^*-1}$, so $V(x)u^{k-1}$ can be seen as an `intermediate' term between the pure reaction term and the pure critical term. In particular, the conclusion of Theorem \ref{mainthm} holds in the following relevant cases:
\begin{enumerate}[label={(\alph*)}]
\item the purely critical problem
$$ -\Delta_p u = \Lambda u^{p^*-1} \quad \mbox{in} \;\;\R^N, $$
\item the weighted eigenvalue problem
$$ -\Delta_p u = V(x)u^{p-1} \quad \mbox{in}\;\;\R^N. $$
\end{enumerate}
\end{rmk}

The aim of the present paper is twofold:
\begin{itemize}
\item generalizing the results in \cite{V} by considering differential inequalities and non-autonomous terms, incorporating a large class of reactions (see Remark \ref{refteo} above);
\item furnishing a quantitative decay estimate for solutions to general singular, convective, critical problems, as the one considered in \cite{BG}, whose prototype is
$$\left\{
\begin{alignedat}{2}
-\Delta_p u &=\lambda (1+|x|)^{-l}(u^{-\gamma}+|\nabla u|^{r-1}) +  u^{p^*-1} \quad&&\mbox{in} \;\; \R^N\\
u &> 0 \quad &&\mbox{in} \;\; \R^N, \\
u(x) &\to 0 \quad &&\mbox{as} \;\; |x|\to+\infty,
\end{alignedat}
\right.$$
where $0<\gamma<1<r<p<N$, $l> N+\gamma \, \frac{N-p}{p-1}$, and $\lambda>0$ is a small parameter.
\end{itemize}

The analysis of the decay rate of solutions to critical problems has a long history. Up to our knowledge, the first contribution in this regard was given by Jannelli and Solimini \cite{JS}, which proved an upper bound estimate for solutions to critical Laplacian problems. After more than fifteen years, Vétois \cite{V} furnished both upper and lower estimates for solutions of a larger class of critical equations governed by the $p$-Laplacian operator. The significant contribution of Vétois paved the way to generalizations for different type of operators and/or nonlinearities: for instance, see \cite{X15, OSV, ELS} for problems exhibiting the Hardy potential, \cite{Vadv} for the anisotropic case, and \cite{BMS16} for the fractional case. Of course, our bibliography is far from being exhaustive.

Theorem \ref{mainthm} provides not only the decay rate of solutions, but also a double-sided decay estimate on their gradients, generalizing Theorem 1.1 in \cite{V}. Moving from the pioneering paper \cite{gnn}, concerning linear problems (see \cite{Li, DR, dpr, sz} for $1<p\leq 2$), the additional information on the gradient decay rate has been successfully combined with the moving planes method to prove radial symmetry results for $C^1$ solutions of autonomous $p$-superlinear problems of the form
\begin{equation*}
\left\{
\begin{alignedat}{2}
-\Delta_p u &=f(u) \quad&&\mbox{in} \;\; \R^N\\
u &> 0 \quad &&\mbox{in} \;\; \R^N, \\
u(x) &\to 0 \quad &&\mbox{as} \;\; |x|\to+\infty.
\end{alignedat}
\right.
\end{equation*}
As customary, symmetry results precede classification results. In the case of the pure critical equation
$$-\Delta_p u=u^{p^*-1}\quad \text{in}\,\, \R^N, \quad u>0\quad \text{in}\,\, \R^N,$$
the classification was given in \cite{cgs} in the linear case $p=2$, while the case $1<p<2$, partially studied by \cite{DMMS}, was completely understood in \cite{V}. Finally, \cite{S16} covered the full range $1<p<N$ (see also \cite{GV}).

Let us sketch the proof of Theorem \ref{mainthm}. Following \cite{V}, our approach is based on a priori estimates, regularity results, and scaling arguments. Roughly speaking, a priori estimates in $L^\infty(\R^N)$ (Lemma \ref{boundlemma}) and $C^{1,\tau}_\loc(\R^N)$ (Lemma \ref{reglemma})
ensure boundedness and regularity of solutions to \eqref{prob} (Theorem \ref{regthm}) and allow to infer a preliminary decay estimate (Theorem \ref{prelimdecay}) of type
$$ u(x) \lesssim |x|^{\frac{p-N}{p}} \quad \mbox{as} \;\; |x|\to+\infty, $$
obtained by using a suitable blow-up sequence given by the doubling lemma (Lemma \ref{lem_dou}). Then an a priori estimate in the weak Lebesgue space $L^{p_*-1,\infty}(\R^N)$ (Lemma \ref{weaklebboundlemma}), being $p_*:=\frac{(N-1)p}{p-1}$, and a scaling argument ensure the decay estimates in Theorem \ref{mainthm}. We also provide the specific dependence of the constants $C_1,C_2$ on the solution $u$, showing that our estimates are uniform when no concentration of compactness (neither at points nor at infinity) occurs (Remark \ref{unifests}).

Accordingly, a key ingredient in our proof is the doubling lemma. Stated for the first time in \cite{hu96} to estimate the blow-up rate of local solutions to parabolic problems, it has been adapted to the elliptic setting in \cite{Souplet}, which provides a universal (i.e., independent of both $u$ and boundary conditions) a priori estimate for non-negative local classical solutions $u\in C^1(\Omega)$ to the Lane-Emden equation
$$ -\Delta_p u = u^q \quad \mbox{in} \;\; \Omega, \quad p-1<q<p^*, $$
that is,
\begin{equation}\label{univest1}
u(x)+|\nabla u(x)|^{\frac{p}{q+1}}\lesssim \dist(x,\partial\Omega)^{-\frac{p}{q+1-p}} \quad \mbox{for all} \;\; x\in\Omega.
\end{equation}
Here $\Omega\subseteq\R^N$ is an arbitrary domain (with $\Omega\neq\R^N$), while $\dist$ represents the distance function. Actually, \cite{Souplet} considers more general nonlinearities, also exhibiting convection terms, but producing, in the general case, the weaker estimate
\begin{equation}
\label{univest2}
u(x)+|\nabla u(x)|^{\frac{p}{q+1}}\lesssim 1+\dist(x,\partial\Omega)^{-\frac{p}{q+1-p}} \quad \mbox{for all} \;\; x\in\Omega.
\end{equation}
Decay estimates in the spirit of \eqref{univest2} have been established for a broader class of operators and nonlinearities: see, e.g., \cite{bfdcds,FSZ}, where also convective terms appear in the reaction term. By the way, such estimates may be obtained by using other techniques: for instance, the direct Bernstein method and a Harnack-type theorem have been exploited in \cite{BV_HG_V} and \cite{szActa}, respectively. 
Note that the universal a priori estimate \eqref{univest2} differs from the uniform $L^\infty$ estimate given by the well-known blow-up method by Gidas and Spruck \cite{gisp}. This latter type of estimates has been generalized in many other contexts, such as for the $N$-Laplacian or for higher-order semilinear problems (see \cite{Rom, ManRom}). By the way, both universal estimates and $L^\infty$ estimates are used, together with fixed point theorems, to get existence results (see \cite{Ruiz,BDGQ,BFNA}).

From a technical point of view, our arguments differ from the ones in \cite{V} mainly for two reasons: (i) we deal with differential inequalities instead of equalities; (ii) the presence of the non-autonomous term $f(x)$ breaks the scaling of the inequality (see Remark \ref{remtec}). In addition, the quantitative estimates furnished in Section 3 allow us to infer uniform decay estimates in absence of concentration phenomena (see Remark \ref{unifests}). We emphasize that the presence of the non-autonomous term in \eqref{prob} compels to a more refined analysis, for instance in the regularity argument, which we treat in the more general case $f\in L^q(\R^N)$, $q\in(N,+\infty]$. However, its impact is most significant in the decay estimates, where additional possibilities, hidden in the purely critical equation, now arise (see Remark \ref{remtec}).

The paper is organized as follows. In Section \ref{prel} we recall some basic results, such as the doubling lemma and the local boundedness theorem. Then Section \ref{regsec} contains a priori estimates and regularity results for solutions to non-autonomous critical problems as \eqref{prob}. Finally, Section \ref{dec} is devoted to the proof of Theorem \ref{mainthm}, chiefly obtained as improvement of a preliminary decay estimate.

\section{Preliminaries}
\label{prel}

\subsection{Notation}

Given any $r\in(1,\infty]$, we denote by $r':=\frac{r}{r-1}$ the Young conjugate of $r$ (with the position $r'=1$ if $r=\infty$). Moreover, $p^*:=\frac{Np}{N-p}$ indicates the Sobolev conjugate of $p$. We will also use $p_*:=\frac{(N-1)p}{N-p}$.

We indicate with $B_r(x)$ the $\R^N$-ball of center $x\in\R^N$ and radius $r>0$, omitting $x$ when it is the origin. The symbols $\overline{B}$, $\partial B$, $B^e$ stand, respectively, for the closure, the boundary, and the exterior of the ball $B$. Given any $K,\Omega\subseteq\R^N$, we write $K\Subset\Omega$ to signify that the closure of $K$ is contained in $\Omega$. For any $N$-dimensional Lebesgue measurable set $\Omega$, by $|\Omega|$ we mean its $N$-dimensional Lebesgue measure. Moreover, for any $x\in \R^N$ we denote with $d(x,\Omega):=\inf\{d(x,y): \, y\in\Omega\}$ the Euclidean distance of $x$ from a set $\Omega$.

Given a real-valued function $\varphi$, we indicate its positive (resp., negative) part with $\varphi_+:=\max\{\varphi,0\}$ (resp., $\varphi_-:=\max\{-\varphi,0\}$). We abbreviate with $\{u>v\}$ the set $\{x\in\R^N: \, u(x)>v(x)\}$, and similarly for $\{u<v\}$, etc.

If a sequence $\{u_n\}$ strongly converges to $u$, we write $u_n\to u$; if the convergence is in weak sense, we use $u_n\rightharpoonup u$. The pair $\langle\cdot,\cdot\rangle$ denotes the duality brackets.

Given any measurable set $\Omega\subseteq \R^N$ and $q\in[1,\infty]$, $L^q(\Omega)$ stands for the standard Lebesgue space, whose norm will be indicated with $\|\cdot\|_{L^q(\Omega)}$, or simply $\|\cdot\|_q$ when $\Omega=\R^N$.

For any $s\in(0,\infty)$, we define the weak Lebesgue space $L^{s,\infty}(\Omega)$ as the set of all measurable functions $u:\Omega\to\R$ such that
$$\|u\|_{L^{s,\infty}(\Omega)}:=\sup_{h>0}\left(h \left|\left\{\left|u\right|>h\right\}\right|^{1/s}\right)<\infty.$$
The quantity $\|\cdot\|_{L^{s,\infty}(\Omega)}$, indicated simply as $\|\cdot\|_{s,\infty}$ when $\Omega=\R^N$, makes $L^{s,\infty}(\Omega)$ a quasi-normed space (see for instance \cite{Gra}). Moreover, the embedding $L^s(\Omega)\hookrightarrow L^{s,\infty}(\Omega)$ is continuous (see \cite[Proposition 1.1.6]{Gra}). The same holds for
\begin{equation}
\label{embedding}
L^{s,\infty}(\Omega)\hookrightarrow L^{s-\eps}(\Omega) \quad \mbox{for all} \;\; \eps\in(0,s),
\end{equation}
provided $\Omega$ has finite measure (see, e.g., \cite[Exercise 1.1.11]{Gra}). Incidentally, we recall the following interpolation inequality (see \cite[Proposition 1.1.14]{Gra}): given any $0<p,q\leq\infty$,
\begin{equation}
\label{interpolation}
\|f\|_r \leq C \|f\|_{p,\infty}^\theta \|f\|_{q,\infty}^{1-\theta} \quad \mbox{for all} \;\; f\in L^{p,\infty}(\R^N)\cap L^{q,\infty}(\R^N),
\end{equation}
where $\frac{1}{r}=\frac{\theta}{p}+\frac{1-\theta}{q}$ and $C>0$ is a suitable constant depending on $p,q,r$.

Let $p\in(1,N)$. Aside from the standard Sobolev spaces $W^{1,p}(\Omega)$, $\Omega\subseteq\R^N$, we will also make use of the Beppo Levi space $\D^{1,p}_0(\R^N)$, which is the closure of the set $C^\infty_c(\R^N)$ of the compactly supported test functions with respect to the norm $\|u\|_{\D^{1,p}_0(\R^N)}:=\|\nabla u\|_p$.
Sobolev's theorem ensures that $\D^{1,p}_0(\R^N) \hookrightarrow L^{p^*}(\R^N)$; the best constant $c$ in the Sobolev inequality $\|u\|_{L^{p^*}(\R^N)}\leq c \|u\|_{\D^{1,p}_0(\R^N)}$ is $S^{-1/p}$, being
$$ S:= \inf_{u\in \D^{1,p}_0(\R^N)\setminus\{0\}} \frac{\|\nabla u\|_p^p}{\|u\|_{p^*}^p}. $$
According to Sobolev's theorem, one has
$\D^{1,p}_0(\R^N) = \left\{u\in L^{p^*}(\R^N): \, |\nabla u|\in L^p(\R^N)\right\}$.

Finally, given any $\Omega\subseteq\R^N$, we indicate with $C^0(\Omega)$ (resp., $C^1(\Omega)$) the set of all continuous (resp., continuously differentiable) functions in $\Omega$. The same holds for $C^{0,\tau}(\Omega)$ and $C^{1,\tau}(\Omega)$, where the continuity is replaced by the H\"older-continuity of exponent $\tau$. We endow those spaces with the norm
\begin{equation*}
\begin{alignedat}{2}
&\|u\|_{C^0(\Omega)} := \sup_\Omega |u|, \quad &&\|u\|_{C^1(\Omega)} := \sup_\Omega |u| + \sup_\Omega |\nabla u|, \\
&\|u\|_{C^{0,\tau}(\Omega)} := \|u\|_{C^0(\Omega)}+[u]_{C^{0,\tau}(\Omega)}, \quad &&\|u\|_{C^{1,\tau}(\Omega)} := \|u\|_{C^1(\Omega)}+[\nabla u]_{C^{0,\tau}(\Omega)},
\end{alignedat}
\end{equation*}
being $[f]_{C^{0,\tau}(\Omega)}$ the semi-norm
$$ [f]_{C^{0,\tau}(\Omega)}:=\sup_{\stackrel{x,y\in\Omega}{x\neq y}} \frac{|f(x)-f(y)|}{|x-y|^\tau}. $$

For every $\Omega\subseteq\R^N$ and $X(\Omega)$ real-valued function space on $\Omega$, by $X(\Omega)_+$ we mean the subspace of $X(\Omega)$ consisting of its non-negative functions. Moreover, we write $u\in X_\loc(\Omega)$ if $u_{\mid_K}\in X(K)$ for all compact $K\Subset\Omega$.

In the sequel, the letter $C$ denotes a positive constant which may change its value at each passage; subscripts on $C$ emphasize its dependence on the specified parameters.

\subsection{Tools}

The following doubling lemma will be crucial in our blow-up arguments.

\begin{lemma}[{\cite[Theorem 5.1]{Souplet}}]\label{lem_dou}
Let $(X,d)$ be a complete metric space, and let $\emptyset \neq D \subset \Sigma \subset X$, with $\Sigma$ closed. Set $\Gamma:=\Sigma \setminus D$. Finally, let $k>0$ and $M: D \to (0,+\infty)$ be bounded on compact subsets of $D$. If $x\in D$ is such that
$$M(x)>\frac{2k}{d(x,\Gamma)},$$
then there exists $y\in D$ such that
\begin{equation}\label{cond_dou}
M(y)>\frac{2k}{d(y,\Gamma)}, \quad M(y)\ge M(x),
\end{equation}
and
$$M(z)\le 2M(y) \quad \mbox{for all} \;\; z\in D\cap \overline{B}_X(y,kM^{-1}(y)),$$
where $\overline{B}_X(x_0,r)$ is the closure of the $X$-ball having center $x_0$ and radius $r$.
\end{lemma}

\begin{rmk}
Let $X=\R^N$ and $\Omega$ be an open subset of $\R^N$. Put $D:=\Omega$ and $\Sigma:=\overline{\Omega}$, so that $\Gamma=\partial\Omega$. Then $\overline{B}_X(y,kM^{-1}(y))\subset D$ for all $y\in D$. Indeed, since $D$ is open, $\eqref{cond_dou}$ implies that
$$d(y,X\setminus D)=d(y,\Gamma)>2kM^{-1}(y).$$
\end{rmk}

We will also use the following version of the well-known local boundedness theorem \cite[Theorem 7.1.1]{PS}.
\begin{thm}
\label{locboundthm}
Let $\Omega\subseteq\R^N$ be a bounded domain and $u\in W^{1,p}_\loc(\R^N)$ be a non-negative solution to
$$ -\Delta_p u \leq au^{p-1}+b \quad \mbox{in} \;\; \Omega. $$
Then, for all $s>p-1$, $x\in\Omega$, and $r>0$ such that $B_{2r}(x)\subseteq\Omega$, one has
$$ \sup_{B_r(x)} u \leq C\left[ r^{-\frac{N}{s}}\|u\|_{L^s(B_{2r}(x))}+r^{p'}b \right], $$
where $C>0$ depends on $p,N,r,s,a$.
\end{thm}

We conclude this section by giving a lower bound for the Beppo Levi norm of non-negative solutions to critical differential inequalities.

\begin{lemma}
\label{univbound}
Let $v\in\D^{1,p}_0(\R^N)\setminus\{0\}$ be a non-negative solution to
$$-\Delta_p v \leq\Lambda  v^{p^*-1} \quad \mbox{in} \;\; \R^N.$$
Then
$$\|\nabla v\|_p \geq \Lambda^{-(N-p)/p^2} S^{N/p^2}=:\kappa_0>0.$$
\end{lemma}
\begin{proof}
Testing with $v$ and using Sobolev's embedding one has
$$ \|\nabla v\|_p^p \leq\Lambda \|v\|_{p^*}^{p^*} \leq \Lambda S^{-p^*/p} \|\nabla v\|_p^{p^*}. $$
Then dividing by $\Lambda S^{-p^*/p}\|\nabla v\|_p^p>0$ yields the conclusion.
\end{proof}

\section{A priori estimates}
\label{regsec}

In this section we recall some classical results concerning essential boundedness, $C^{1,\tau}$ local regularity, and $p^*$-weak summability of solutions to $p$-Laplacian problems, here adapted to the setting of \eqref{prob}.

\begin{lemma}\label{boundlemma}
Let $f\in L^q(\R^N)_+$ for some $q\in\left(\frac{N}{p},\infty\right]$. Then any solution $u\in\D^{1,p}_0(\R^N)$ of \eqref{prob} belongs to $L^\infty(\R^N)$.
\end{lemma}
\begin{proof}
Set $\eta:=\left(\frac{p}{N}-\frac{1}{q}\right)^{-1}>0$, take any $k>0$ and let $\Omega_k:=\{u>k\}$. Since $u\in L^{p^*}(\R^N)$, then the Chebichev inequality entails 
\begin{equation}
\label{vanishsuperlevel}
|\Omega_k| \leq k^{-p^*} \|u\|_{p^*}^{p^*} \to 0 \quad \mbox{as} \;\; k\to+\infty.
\end{equation}
Hence, there exists $\overline{k}>1$ such that $|\Omega_k|^\eta<\frac{S}{2\|f\|_q}$ for all $k>\overline{k}$, being $S>0$ the Sobolev constant. Now pick any $k>\overline{k}$. Testing \eqref{prob} with $(u-k)_+$, besides exploiting $f\in L^q(\R^N)$, as well as H\"older's and Sobolev's inequalities, we get
\begin{equation*}
\begin{aligned}
\|\nabla (u-k)\|_{L^p(\Omega_k)}^p &\leq \int_{\Omega_k} fu \dx + \Lambda\|u\|_{L^{p^*}(\Omega_k)}^{p^*} \leq \int_{\Omega_k} fu^p \dx + \Lambda\|u\|_{L^{p^*}(\Omega_k)}^{p^*} \\
&\leq \|f\|_q \|u\|_{L^{p^*}(\Omega_k)}^p|\Omega_k|^\eta + \Lambda\|u\|_{L^{p^*}(\Omega_k)}^{p^*} \\
&\leq \frac{\|f\|_q}{S} \, \|\nabla(u-k)\|_{L^p(\Omega_k)}^p |\Omega_k|^\eta + \Lambda\|u\|_{L^{p^*}(\Omega_k)}^{p^*} \\
&\leq \frac{1}{2} \, \|\nabla(u-k)\|_{L^p(\Omega_k)}^p + 2^{p^*-1}\Lambda\left(\|u-k\|_{L^{p^*}(\Omega_k)}^{p^*}+k^{p^*}|\Omega_k|\right).
\end{aligned}
\end{equation*}
Reabsorbing $\frac{1}{2} \|\nabla(u-k)\|_{L^p(\Omega_k)}^p$ on the left and using Sobolev's inequality again, we obtain
$$ \|u-k\|_{L^{p^*}(\Omega_k)}^p \leq C\left(\|u-k\|_{L^{p^*}(\Omega_k)}^{p^*}+k^{p^*}|\Omega_k|\right), $$
for some $C>0$ depending only on $p,N,q,\|f\|_q,\Lambda,u$.

Let $M>2\overline{k}$ and set $k_n:=M(1-2^{-n})$ for all $n\in\N$. Repeating verbatim the proof of \cite[Lemma 3.2]{CGL}, we infer that 
\begin{equation}\label{kn}
\|u-k_n\|_{L^{p^*}(\Omega_{k_n})} \to 0 \quad \mbox{as} \;\; n\to +\infty,
\end{equation}
provided
\begin{equation}
\label{equiuniformint}
\|u-M/2\|_{L^{p^*}(\Omega_{M/2})}
\end{equation}
is small enough (that can be ensured by taking $M$ large enough, depending on $u$). By the definition of $k_n$ and \eqref{kn}, one has
\begin{equation*}
\int_{\R^N} \left(u-M\right)_+^{p^*} \dx\le \int_{\R^N} \left(u-k_n\right)_+^{p^*} \dx=\int_{\Omega_{k_n}} \left(u-k_n\right)^{p^*} \dx  \to 0 \quad \mbox{as} \;\; n\to +\infty.
\end{equation*}
Hence $u\leq M$ almost everywhere in $\R^N$, so that $u\in L^\infty(\R^N)$.
\end{proof}

\begin{lemma}
\label{reglemma}
Let $\Omega\subseteq \R^N$ be an open set, $q\in(N,\infty]$, and $u\in W^{1,p}_\loc(\Omega)$ be a solution to
\begin{equation}
\label{doubleineq}
0\leq -\Delta_p u\leq h(x) \quad \mbox{in} \;\; \Omega
\end{equation}
for some $h\in L^q_\loc(\Omega)_+$. Then $u\in C^{1,\tau}_\loc(\Omega)$ for some $\tau\in(0,1)$. Moreover,
\begin{equation}
\label{regest}
\|\nabla u\|_{C^{0,\tau}(B_R(x))} \leq C\left(\|\nabla u\|_{L^p(B_{2R}(x))}^{p-1}+\|h\|_{L^q(B_{2R}(x))}\right)
\end{equation}
for all $B_{2R}(x)\Subset \Omega$, where $C$ depends on $p,N,q,R$.
\end{lemma}

\begin{proof}
Let $x\in\Omega$ and $R>0$ such that $B_{2R}(x)\Subset\Omega$. Hereafter, we will omit the dependence on $x$. Observe that $C^\infty_c(B_{2R})\hookrightarrow L^{q'}(B_{2R})$ densely (see \cite[Corollary 4.23]{B}), whence
\begin{equation}
\label{Lqembedding}
\|\phi\|_{L^{q'}(B_{2R})} \leq C_R \|\phi\|_{C^\infty_c(B_{2R})} \quad \mbox{for all} \;\; \phi\in C^\infty_c(B_{2R}),
\end{equation}
for a suitable $C_R>0$ depending on $N,R,q$. Let us consider $-\Delta_p u \in (C^\infty_c(B_{2R}))^*$ and take any $\phi\in C^\infty_c(B_{2R})$. Testing \eqref{doubleineq} with $\phi_\pm$ gives
$$ 0 \leq \langle -\Delta_p u,\phi_\pm \rangle \leq \int_{B_{2R}} h\phi_\pm \dx. $$
By linearity, besides \eqref{Lqembedding} and H\"older's inequality, we obtain
\begin{equation}
\label{normest}
\begin{aligned}
\frac{\left|\langle -\Delta_p u,\phi \rangle\right|}{\|\phi\|_{C^\infty_c(B_{2R})}} &\leq C_R\,\frac{\left|\langle -\Delta_p u,\phi_+ \rangle - \langle -\Delta_p u,\phi_- \rangle\right|}{\|\phi\|_{L^{q'}(B_{2R})}} \\
&\leq C_R\,\frac{\int_{B_{2R}} h\phi_+ \dx + \int_{B_{2R}} h\phi_- \dx}{\|\phi\|_{L^{q'}(B_{2R})}} \leq C_R\|h\|_{L^q(B_{2R})}.
\end{aligned}
\end{equation}
Therefore, according to the Hahn-Banach theorem \cite[Corollary 1.2]{B}, there exists a unique $A_p(u)\in (L^{q'}(B_{2R}))^*$ extending $-\Delta_p u$ and such that
\begin{equation}
\label{hahnbanach}
\|A_p(u)\|_{(L^{q'}(B_{2R}))^*} = \|-\Delta_p u\|_{(C^\infty_c(B_{2R}))^*}.
\end{equation}
Thus, owing to the Riesz theorem \cite[Theorems 4.11 and 4.14]{B}, there exists $g\in L^q(B_{2R})$ such that $A_p(u)=g$, that is,
$$ \langle A_p(u),\psi \rangle = \int_{B_{2R}} g\psi \dx \quad \mbox{for all} \;\; \psi\in L^{q'}(B_{2R}). $$
Thus, recalling \eqref{hahnbanach} and \eqref{normest},
\begin{equation}
\label{normcompare}
\|g\|_{L^q(B_{2R})} = \|A_p(u)\|_{(L^{q'}(B_{2R}))^*} = \|-\Delta_p u\|_{\left(C^\infty_c(B_{2R})\right)^*} \leq C_R \|h\|_{L^q(B_{2R})}.
\end{equation}

According to \cite[Theorem 1.5]{DM} (see also \cite[Lemma 2.4]{GM}) and the inequality
$$\|h\|_{L^2(B_\rho(y))}\leq (|B_1|\rho^N)^{\frac{1}{2}-\frac{1}{q}}\|h\|_{L^q(B_\rho(y))} \quad \mbox{for all} \;\; y\in\R^N \;\; \mbox{and} \;\; \rho>0,$$
for any $z\in B_{\frac{2}{3}R}$ one has
\begin{equation*}
\begin{aligned}
\|\nabla u\|_{L^\infty(B_{R/3}(z))}^{p-1} &\leq C\left( \|\nabla u\|_{L^p(B_{\frac{2}{3}R}(z))}^{p-1} + \sup_{y\in B_{\frac{2}{3}R}(z)} \int_0^{\frac{2}{3}R} \rho^{-\frac{N}{2}}\|h\|_{L^2(B_\rho(y))} \drho \right) \\
&\leq C\left( \|\nabla u\|_{{L^p(B_{\frac{2}{3}R}(z))}}^{p-1} + \|h\|_{L^q(B_{\frac{4}{3}R}(z))} \int_0^{2R} \rho^{-\frac{N}{q}} \drho \right) \\
&\leq C\left( \|\nabla u\|_{L^p(B_{2R})}^{p-1} + \|h\|_{L^q(B_{2R})} \right),
\end{aligned}
\end{equation*}
where $C>0$ depends on $p,N,R,q$. Arbitrariness of $z\in B_{\frac{2}{3}R}$ implies
\begin{equation}
\label{supgradest}
\|\nabla u\|_{L^\infty(B_R)}^{p-1} \leq C\left( \|\nabla u\|_{L^p(B_{2R})}^{p-1} + \|h\|_{L^q(B_{2R})} \right).
\end{equation}

Let us consider $v\in W^{1,2}_0(B_{\frac{3}{2}R})$ the unique solution to
\begin{equation*}
\left\{
\begin{alignedat}{2}
-\Delta v &= g(x) \quad &&\mbox{in} \;\; B_{\frac{3}{2}R}, \\
v &= 0 \quad &&\mbox{on} \;\; \partial B_{\frac{3}{2}R}.
\end{alignedat}
\right.
\end{equation*}
By using Morrey's inequality \cite[Theorem 9.12]{B}, the Calderon-Zygmund theorem \cite[Theorem 9.9]{GT}, and \eqref{normcompare} we infer
\begin{equation}
\label{oscillation}
\|\nabla v\|_{C^{0,\beta}(B_{\frac{3}{2}R})} \leq C\|v\|_{W^{2,q}(B_{\frac{3}{2}R})} \leq C\|g\|_{L^q(B_{\frac{3}{2}R})} \leq C\|h\|_{L^q(B_{2R})},
\end{equation}
where $\beta=1-\frac{N}{q}\in (0,1)$ and $C>0$ depends on $p,N,R,q$. Notice that $-\Delta_p u = \Div(\nabla v)$ in $B_{\frac{3}{2}R}$. Then \cite[Corollary 5.2]{BCDKS} and \eqref{oscillation} ensure
\begin{equation}
\label{cianchiest}
\begin{aligned}
[|\nabla u|^{p-2}\nabla u]_{C^{0,\beta}(B_R)} &\leq C\left(\|\nabla u\|_{L^p(B_{\frac{3}{2}R})}^{p-1} + \|\nabla v\|_{C^{0,\beta}(B_{\frac{3}{2}R})}\right) \\
&\leq C\left(\|\nabla u\|_{L^p(B_{2R})}^{p-1} + \|h\|_{L^q(B_{2R})}\right),
\end{aligned}
\end{equation}
where $C>0$ depends on $p,N,R,q$.

Set $\tau:=\beta\min\{\frac{1}{p-1},1\}$. Let us prove that, for some $C>0$ depending on $p,N,R,q$,
\begin{equation}
\label{oscest}
[\nabla u]_{C^{0,\tau}(B_R)}^{p-1} \leq C\left(\|\nabla u\|_{L^p(B_{2R})}^{p-1} + \|h\|_{L^q(B_{2R})}\right).
\end{equation}
To this aim, we consider the bijective map $ \Psi: \R^N \to \R^N $ defined as
\begin{equation*}
\Psi(y) = \left\{
\begin{array}{ll}
|y|^{\frac{2-p}{p-1}} \, y \quad &\mbox{if} \; \; y \neq 0, \\
0 \quad &\mbox{if} \; \; y = 0,
\end{array}
\right.
\end{equation*}
whose inverse is $ \Psi^{-1}(y) = |y|^{p-2} y $. Inequalities (I) and (VII) of \cite[Chapter 12]{Lind} ensure that, for any $ y_1,y_2 \in \R^N $,
\begin{equation}
\label{lind}
\begin{split}
&(\Psi^{-1}(y_1)-\Psi^{-1}(y_2)) \cdot (y_1-y_2) \\
&\geq \left\{
\begin{array}{ll}
2^{2-p}|y_1-y_2|^p \quad &\mbox{if} \; \; p \geq 2, \\
(p-1)(1+|y_1|^2+|y_2|^2)^{\frac{p-2}{2}} |y_1-y_2|^2 \quad &\mbox{if} \; \; 1<p<2.
\end{array}
\right.
\end{split}
\end{equation}
Let us discuss two cases: $p\geq 2$ and $p\in (1,2)$. \\
{\bf Case I:} \, $p\geq 2$. Pick any $x_1,x_2\in B_R$ with $x_1\neq x_2$. Applying the Cauchy-Schwartz inequality to \eqref{lind} with $y_1:=\nabla u(x_1)$ and $y_2:=\nabla u(x_2)$ one has
$$ \left|\nabla u(x_1)-\nabla u(x_2))\right|^{p-1} \leq 2^{p-2} \left||\nabla u(x_1)|^{p-2}\nabla u(x_1)-|\nabla u(x_2)|^{p-2}\nabla u(x_2)\right|. $$
Therefore, exploiting also \eqref{cianchiest},
$$ [\nabla u]_{C^{0,\tau}(B_R)}^{p-1} \leq 2^{p-2} [|\nabla u|^{p-2}\nabla u]_{C^{0,\beta}(B_R)} \leq C\left(\|\nabla u\|_{L^p(B_{2R})}^{p-1} + \|h\|_{L^q(B_{2R})}\right). $$
{\bf Case II:} \, $p\in(1,2)$. It is not restrictive to assume $ \|h\|_{L^q(B_{2R})} = 1 $: indeed, assume that \eqref{oscest} holds when the datum of \eqref{doubleineq} has unitary $L^q(B_{2R})$-norm, and set $\lambda := \|h\|_{L^q(B_{2R})}^{\frac{1}{1-p}}$, so that
\begin{equation*}
-\Delta_p (\lambda u) = \lambda^{p-1} (-\Delta_p u) \leq \lambda^{p-1} h
\end{equation*}
implies
\begin{equation*}
[\nabla (\lambda u)]_{C^{0,\tau}(B_R)}^{p-1} \leq C\left(\|\nabla (\lambda u)\|_{L^p(B_{2R})}^{p-1} + \|\lambda^{p-1} h\|_{L^q(B_{2R})}\right),
\end{equation*}
whence \eqref{oscest}, dividing by $\lambda^{p-1}$.

Hence, assume $\|h\|_{L^q(B_{2R})}=1$. Reasoning as above, for any $x_1,x_2\in B_R$, \eqref{lind} gives
\begin{equation*}
\begin{aligned}
&|\nabla u(x_1) - \nabla u(x_2)| \\
&\leq \frac{1}{p-1}(1+|\nabla u(x_1)|^2+|\nabla u(x_2)|^2)^{\frac{2-p}{2}}\left||\nabla u(x_1)|^{p-2}\nabla u(x_1)-|\nabla u(x_2)|^{p-2}\nabla u(x_2)\right| \\
&\leq \frac{1}{p-1}\left(1+2\|\nabla u\|_{L^\infty(B_R)}^2\right)^{\frac{2-p}{2}}\left||\nabla u(x_1)|^{p-2}\nabla u(x_1)-|\nabla u(x_2)|^{p-2}\nabla u(x_2)\right|,
\end{aligned}
\end{equation*}
where we have used \eqref{supgradest} to ensure that $\|\nabla u\|_{L^\infty(B_R)}$ is finite. Then, setting
$$ \Theta:= \|\nabla u\|_{L^p(B_{2R})}^{p-1}+\|h\|_{L^q(B_{2R})}\geq \|h\|_{L^q(B_{2R})} = 1 $$
as well as exploiting \eqref{supgradest} and \eqref{cianchiest}, we get
\begin{equation*}
\begin{split}
[\nabla u]_{C^{0,\tau}(B_R)} &\leq \frac{1}{p-1} \left(1+2\|\nabla u\|_{L^\infty(B_R)}^2\right)^{\frac{2-p}{2}} [|\nabla u|^{p-2} \nabla u]_{C^{0,\tau}(B_R)} \\
&\leq \frac{C}{p-1} \left[1+C \Theta^\frac{2}{p-1}\right]^{\frac{2-p}{2}} \Theta \\&\leq C\Theta^{\frac{1}{p-1}} = C\left(\|\nabla u\|_{L^p(B_{2R})}^{p-1} + \|h\|_{L^q(B_{2R})}\right)^{\frac{1}{p-1}},
\end{split}
\end{equation*}
yielding \eqref{oscest}.

Putting \eqref{supgradest} and \eqref{oscest} together entail \eqref{regest}, concluding the proof.
\end{proof}

\begin{thm}\label{regthm}
Let $u\in\D^{1,p}_0(\R^N)$ be a solution to \eqref{prob} with $f\in L^q(\R^N)_+$ for some $q\in(N,\infty]$. Then $u\in L^\infty(\R^N)\cap C^{1,\tau}_\loc(\R^N)$ for some $\tau\in(0,1)$.
\end{thm}
\begin{proof}
This is a direct consequence of Lemmas \ref{boundlemma}--\ref{reglemma}. Indeed, Sobolev's inequality and Lemma \ref{boundlemma} ensure $u\in L^{p^*}(\R^N)\cap L^\infty(\R^N)$, while \eqref{interpolation} and $(p^*)'<N<q\leq \infty$ yield
$$ \left\|f(x)+\Lambda u^{p^*-1}\right\|_q \leq \|f\|_q+\Lambda \|u\|_{p^*}^{\frac{p^*}{q}}\|u\|_\infty^{\frac{p^*}{q'}-1}<+\infty, $$
leading to the conclusion via Lemma \ref{reglemma}.
\end{proof}

\begin{lemma}
\label{weaklebboundlemma}
Let $f\in L^1(\R^N)_+\cap L^q(\R^N)_+$ for some $q\in\left(\frac{N}{p},\infty\right]$. Then any solution $u$ to \eqref{prob} belongs to $L^{p_*-1,\infty}(\R^N)$. In particular, $u\in L^r(\R^N)$ for all $r\in(p_*-1,\infty]$.
\end{lemma}
\begin{proof}
Take any $u$ solving \eqref{prob}. Theorem \ref{regthm} ensures that $u\in L^\infty(\R^N)$, so we put $M:=\|u\|_\infty$. For any $h\in(0,M)$, test \eqref{prob} with $T_h(u):=\min\{u,h\}$. Then
\begin{equation}
\label{Thtest}
\int_{\{u\leq h\}} |\nabla u|^p \dx \leq \int_{\R^N} f T_h(u) \dx + \Lambda\left[\int_{\{u\leq h\}} u^{p^*} \dx + h\int_{\{u>h\}} u^{p^*-1} \dx\right].
\end{equation}
Now we estimate the terms into square brackets. It is readily seen that
\begin{equation}
\label{estterm1}
\int_{\{u\leq h\}} u^{p^*} \dx = \int_{\R^N} T_h(u)^{p^*} \dx - h^{p^*}|\{u>h\}|,
\end{equation}
while the layer-cake lemma \cite[Proposition 1.1.4]{Gra} entails
\begin{equation}
\label{estterm2}
\begin{aligned}
\int_{\{u>h\}} u^{p^*-1} \dx &= (p^*-1)\int_0^M s^{p^*-2} |\{u>\max\{s,h\}\}| \ds \\
&= h^{p^*-1}|\{u>h\}|+(p^*-1)\int_h^M s^{p^*-2}|\{u>s\}| \ds.
\end{aligned}
\end{equation}
Thus, plugging \eqref{estterm1}--\eqref{estterm2} into \eqref{Thtest} yields
\begin{equation}
\label{Thtest2}
\begin{aligned}
&\int_{\{u\leq h\}} |\nabla u|^p \dx \\
&\leq \int_{\R^N} f T_h(u) \dx + \Lambda\left[\int_{\R^N} T_h(u)^{p^*} \dx + (p^*-1)h\int_h^M s^{p^*-2}|\{u>s\}| \ds\right].
\end{aligned}
\end{equation}
In order to reabsorb $\int_{\R^N} f T_h(u) \dx$ into $h\int_h^M s^{p^*-2}|\{u>s\}| \ds$, we notice that
\begin{equation}
\label{layercakelimit}
\lim_{h\to 0^+} \int_h^M s^{p^*-1}|\{u>s\}| \ds = \frac{\|u\|_{p^*}^{p^*}}{p^*}
\end{equation}
by the layer-cake lemma, so there exists $\overline{h}>0$ (depending on $u$) such that 
$$\int_h^M s^{p^*-1}|\{u>s\}| \ds > \frac{\|u\|_{p^*}^{p^*}}{2p^*} \quad \mbox{for all} \;\; h\in(0,\overline{h}).$$
Hence, exploiting also the H\"older inequality,
\begin{equation}
\label{estterm3}
\begin{aligned}
h\int_h^M s^{p^*-2}|\{u>s\}| \ds &\geq \frac{h}{M} \int_h^M s^{p^*-1}|\{u>s\}| \ds > \frac{\|u\|_{p^*}^{p^*}}{2Mp^*}\,h \\
&= \frac{\|u\|_{p^*}^{p^*}}{2Mp^*\|f\|_1}\,\|f\|_1 h \geq \frac{\|u\|_{p^*}^{p^*}}{2Mp^*\|f\|_1}\int_{\R^N} fT_h(u) \dx
\end{aligned}
\end{equation}
for every $h\in(0,\overline{h})$. Putting \eqref{Thtest2} and \eqref{estterm3} together yields
$$ \int_{\{u\leq h\}} |\nabla u|^p \dx \leq C\left[\int_{\R^N} T_h(u)^{p^*} \dx + h\int_h^M s^{p^*-2}|\{u>s\}| \ds\right] $$
for some $C>0$ depending on $p,N,\|f\|_1,\Lambda,u$. Consequently, through Sobolev's inequality we obtain
\begin{equation}
\label{Thtest3}
\begin{aligned}
\int_{\R^N} T_h(u)^{p^*} \dx &\leq S^{-\frac{N}{N-p}} \left(\int_{\{u\leq h\}} |\nabla u|^p \dx\right)^{\frac{N}{N-p}} \\
&\leq C\left[\int_{\R^N} T_h(u)^{p^*} \dx + h\int_h^M s^{p^*-2}|\{u>s\}| \ds\right]^{\frac{N}{N-p}} \\
&\leq C\left[\left(\int_{\R^N} T_h(u)^{p^*} \dx\right)^{\frac{N}{N-p}} + \left(h\int_h^M s^{p^*-2}|\{u>s\}| \ds\right)^{\frac{N}{N-p}}\right]
\end{aligned}
\end{equation}
for all $h\in(0,\overline{h})$. Lebesgue's dominated convergence theorem and $u\in L^{p^*}(\R^N)$ entail $\int_{\R^N} T_h(u)^{p^*} \dx \to 0$ as $h\to 0^+$, so taking a smaller $\overline{h}$ we have
\begin{equation}
\label{reabsorbTh}
C\left(\int_{\R^N} T_h(u)^{p^*} \dx\right)^{\frac{N}{N-p}} \leq \frac{1}{2}\int_{\R^N} T_h(u)^{p^*} \dx \quad \mbox{for all} \;\; h\in(0,\overline{h}).
\end{equation}
Thus, reabsorbing $\int_{\R^N} T_h(u)^{p^*} \dx$ to the left, from \eqref{Thtest3} we infer
\begin{equation*}
\int_{\R^N} T_h(u)^{p^*} \dx \leq C\left(h\int_h^M s^{p^*-2}|\{u>s\}| \ds\right)^{\frac{N}{N-p}} \quad \mbox{for all} \;\; h\in (0,\overline{h}).
\end{equation*}
We deduce
\begin{equation*}
h^{p^*}|\{u>h\}| \leq \int_{\R^N} T_h(u)^{p^*} \dx \leq C\left(h\int_h^M s^{p^*-2}|\{u>s\}| \ds\right)^{\frac{N}{N-p}} \quad \mbox{for all} \;\; h\in (0,\overline{h}).
\end{equation*}
Reasoning as in \cite[p.154]{V}, it turns out that $u\in L^{p_*-1,\infty}(\R^N)$. Then \eqref{interpolation} and $u\in L^\infty(\R^N)$ ensure $u\in L^r(\R^N)$ for all $r\in(p_*-1,\infty]$.
\end{proof}

\section{Decay estimates}
\label{dec}

Before proving Theorem \ref{mainthm}, we preliminarly produce a (non-optimal) decay estimate for solutions to \eqref{prob}.

\begin{thm}\label{prelimdecay}
Let \eqref{hyp} be satisfied and $\kappa_0$ be as in Lemma \ref{univbound}. Then, for any $\kappa<\kappa_0$, $\kappa'>0$, and $r>r'>0$, there exists $K_0>0$ (depending on $\kappa,\kappa',r,r',p,N,\Lambda,L,\alpha$) such that
\begin{equation}
\label{decay}
u(x)\leq K_0|x|^{\frac{p-N}{p}} \quad \mbox{in} \;\; B_r^e
\end{equation}
for every $u\in\D^{1,p}_0(\R^N)$ solution to \eqref{prob} satisfying
$$ \|\nabla u\|_{L^p(B_{r'}^e)} \leq \kappa \quad \mbox{and} \quad \|\nabla u\|_p \leq \kappa'. $$
\end{thm}

\begin{proof}
Fix $\kappa,\kappa',r,r'$ as in the statement, and set $r'':=\frac{r+r'}{2}$. Notice that \eqref{decay} is ensured by
\begin{equation}\label{cladec}
d(x,B_{r''})u(x)^{\frac{p}{N-p}} \leq K' \quad \mbox{for all}\;\; x\in B_r^e,
\end{equation}
where $K'>0$ is a suitable constant depending on $k,k',r,r',p,N,\Lambda,L,\alpha$. Indeed, if \eqref{cladec} holds true, then there exists $K''=K''(r,r')>0$ such that $d(x,B_{r''})>K''|x|$ for all $x\in B_r^e$, since
$$ |x|>r \quad\Rightarrow \quad r''<\frac{r''}{r}\,|x| \quad\Rightarrow \quad d(x,B_{r''})=|x|-r''>\left(1-\frac{r''}{r}\right)|x|=:K''|x|,$$
recalling that $r''<r$. Hence, our aim is to prove \eqref{cladec}.

Assume by contradiction that \eqref{cladec} fails, so that there exist $\{u_n\}\subseteq \D^{1,p}_0(\R^N)$ and $\{x_n\}\subseteq B_r^e$ such that, for all $n\in\N$, $u_n>0$ in $\R^N$ and
\begin{enumerate}[label=${\rm a_{\arabic*})}$, ref=${\rm a_{\arabic*}}$]
\itemsep0.5em
\item \label{equn} $0\leq -\Delta_p u_n \leq f(x)+\Lambda u_n^{p^*-1}$ in $\R^N$,
\item \label{extbound} $\|\nabla u_n\|_{L^p(B_{r'}^e)}\leq \kappa$,
\item \label{sobolevbound} $\|\nabla u_n\|_p \leq \kappa'$,
\item \label{a4} $d(x_n,B_{r''})u_n(x_n)^{\frac{p}{N-p}}>2n$.
\end{enumerate}
According to Theorem \ref{regthm}, $u_n\in L^\infty(\Omega)$ for all $n\in\N$. Then Lemma \ref{lem_dou} (with $k=n$, $M=u_n^{\frac{p}{N-p}}$, and $D= \overline{B}_{r''}^e\subseteq {B}_{r''}^e = \Sigma$) furnishes $\{y_n\}\subseteq  \overline{B}_{r''}^e$ such that
\begin{enumerate}[label=${\rm b_{\arabic*})}$, ref=${\rm b_{\arabic*}}$]
\itemsep0.5em
\item \label{distance} $d(y_n,B_{r''})u_n(y_n)^{\frac{p}{N-p}}>2n$,
\item \label{control} $u_n(x_n)\leq u_n(y_n)$,
\item \label{doubling} $u_n(y)\leq 2^{\frac{N-p}{p}} u_n(y_n)$ for all $y\in B_{n u_n(y_n)^{\frac{p}{p-N}}}(y_n)$.
\end{enumerate}
For every $n\in\N$ we define $\mu_n:=u_n(y_n)^{-1}$ and
\begin{equation}
\label{blowup}
\tilde{u}_n(y) := \mu_n u_n(\mu_n^{\frac{p}{N-p}}y+y_n) \quad \mbox{for all} \;\; y\in\R^N.
\end{equation} 
Thus, from \eqref{equn} we infer
\begin{equation}
\label{equnblow}
0\leq -\Delta_p \tilde{u}_n \leq \mu_n^{p^*-1}f(\mu_n^{\frac{p}{N-p}}y+y_n) + \Lambda\tilde{u}_n^{p^*-1} \quad \mbox{in} \;\; \R^N.
\end{equation}
Exploiting \eqref{doubling}, along with the definition of $\mu_n$ and \eqref{blowup}, yields
\begin{equation}
\label{propblow}
\tilde{u}_n(0)=1 \quad \mbox{and} \quad \tilde{u}_n(y)\leq 2^{\frac{N-p}{p}} \quad \mbox{for all }\;\; y\in B_n.
\end{equation}

We claim that 
\begin{equation}\label{decnonlin}
\lim_{n\to\infty} \mu_n^{p^*-1}f(\mu_n^{\frac{p}{N-p}}y+y_n) = 0 \quad \mbox{locally uniformly in} \;\; y\in\R^N.
\end{equation}
Reasoning up to sub-sequences, we reduce to two cases:
\begin{enumerate}[label={${\rm (\Roman*)}$},ref={${\rm \Roman*}$}]
\itemsep0.5em
\item \label{vanishing} $\mu_n\to 0$;
\item \label{lowerbound} $\mu_n\geq c$ for all $n\in\N$, being $c>0$ opportune.
\end{enumerate}
In case \eqref{vanishing}, it is readily seen that \eqref{decnonlin} holds true uniformly in $\R^N$, since $f\in L^\infty(\R^N)$ by \eqref{hyp}. Otherwise, if \eqref{lowerbound} occurs, fix any compact $H\subseteq\R^N$. Since $B_n\nearrow \R^N$, then $H\subseteq B_n$ for all $n\in\N$ sufficiently large, say $n>\nu\in\N$. Thus, \eqref{distance} and $d(y_n,B_{r''})<|y_n|$ yield, for all $y\in H$,
$$ |\mu_n^{\frac{p}{N-p}}y+y_n| \geq |y_n|-\mu_n^{\frac{p}{N-p}}|y| \geq \mu_n^{\frac{p}{N-p}}(2n-|y|) \geq \mu_n^{\frac{p}{N-p}}n \quad \mbox{for all} \;\; n>\nu. $$
This, together with \eqref{hyp} and $\alpha>N>\frac{N}{p'}+1$, ensures
\begin{equation}
\label{locunifest}
\mu_n^{p^*-1}f(\mu_n^{\frac{p}{N-p}}y+y_n)\le L\mu_n^{-\frac{(\alpha-N)p}{N-p}-1}n^{-\alpha}\le Lc^{-\frac{(\alpha-N)p}{N-p}-1}n^{-\alpha}
\end{equation}
for all $y\in H$ and $n>\nu$. Letting $n\to\infty$, besides taking into account the arbitrariness on $H$, proves \eqref{decnonlin}.

Fix any $x\in\R^N$ and observe that $B_2(x) \subseteq B_n$ for any $n$ sufficiently large. 
Thus, \eqref{propblow} in \eqref{equnblow} gives
\begin{equation}\label{utildeprob}
0\leq -\Delta_p \tilde u_n\leq \mu_n^{p^*-1}f(\mu_n^{\frac{p}{N-p}}y+y_n)+\Lambda2^{\frac{N-p}{p}(p^*-1)} \quad \mbox{in} \;\; B_2(x)
\end{equation}
for all $n$ large enough. Thus \eqref{utildeprob}, together with \eqref{decnonlin}, implies that $\{\tilde u_n\}$ is bounded in $C^{1,\tau}(\overline{B}_1(x))$ for some $\tau\in(0,1)$, according to Lemma \ref{reglemma}. In particular, since $x$ is arbitrary, Ascoli-Arzelà's theorem and a diagonal argument ensure
\begin{equation}
\label{c1locconv}
\tilde u_n\to\tilde{u}_\infty \quad \mbox{in} \;\; C^1_\loc(\R^N)
\end{equation}
for some $\tilde{u}_\infty\in C^1_\loc(\R^N)$. Moreover, \eqref{sobolevbound} gives
\begin{equation}
\label{energybound}
\|\nabla \tilde{u}_n\|_p = \|\nabla u_n\|_p \leq \kappa',
\end{equation}
so that $\tilde{u}_n \rightharpoonup \tilde{u}_\infty$ in $\D^{1,p}_0(\R^N)$.
Then, letting $n\to\infty$ in \eqref{equnblow} and using \eqref{decnonlin}, it turns out that
$$ -\Delta_p \tilde{u}_\infty \leq \Lambda\tilde u_\infty^{p^*-1} \quad \mbox{in} \;\; \R^N. $$
Observe that, for any $R>0$, one has
\begin{equation}\label{localenergy}
\|\nabla \tilde{u}_n\|_{L^p(B_R)} = \|\nabla u_n\|_{L^p(B_{R\mu_n^{\frac{p}{N-p}}}(y_n))} \quad \mbox{for all} \;\; n\in\N. 
\end{equation}

Now we prove that, for all $n\in\N$ large enough,
\begin{equation}\label{disjoint}
B_{R\mu_n^{\frac{p}{N-p}}}(y_n) \cap B_{r'} = \emptyset.
\end{equation}
To this end, take any $y\in B_{r'}$ and observe that \eqref{distance} implies
\begin{equation*}
\begin{aligned}
\mu_n^{-\frac{p}{N-p}}|y-y_n|&\ge \mu_n^{-\frac{p}{N-p}}|y_n|-\mu_n^{-\frac{p}{N-p}}|y| = \mu_n^{-\frac{p}{N-p}}d(y_n,B_{r''})+\mu_n^{-\frac{p}{N-p}}(r''-|y|)\\
&\geq \mu_n^{-\frac{p}{N-p}}d(y_n,B_{r''}) > 2n,
\end{aligned}
\end{equation*}
exploiting also $y_n\in \overline{B}_{r''}^e$ and $y\in B_{r''}$. Therefore, \eqref{disjoint} holds true for all $n>\frac{R}{2}$.

From \eqref{localenergy}--\eqref{disjoint} and \eqref{extbound} we deduce
$$ \|\nabla \tilde{u}_n\|_{L^p(B_R)} \leq \|\nabla u_n\|_{L^p(B_{r'}^e)} \leq \kappa $$
for large $n$'s. Passing to the limit with respect to $n$ through \eqref{c1locconv}, and then letting $R\to+\infty$, we infer
$$ \|\nabla \tilde{u}_\infty\|_p \leq \kappa. $$
Since $\kappa<\kappa_0$, Lemma \ref{univbound} forces $\tilde{u}_\infty=0$, which contradicts $\tilde{u}_\infty(0)=1$, guaranteed by \eqref{propblow} and \eqref{c1locconv}.
\end{proof}

\begin{rmk}
In the proof of Theorem \ref{prelimdecay}, \eqref{hyp} comes into play only to ensure \eqref{decnonlin}, and the requirement $\alpha>N$ can be relaxed to $\alpha\geq \frac{N}{p'}+1$ (cfr. \eqref{locunifest}).
\end{rmk}

\begin{rmk}
The constant $K_0$ depends on $k'$ which, unlike \cite{V}, bounds $\|\nabla u\|_p$ instead $\|u\|_{p^*}$. This is due to the lack of a reverse Sobolev inequality (to be used in \eqref{energybound}), which is guaranteed when $u_n$ fulfills $-\Delta_p u_n \leq\Lambda u_n^{p^*-1}$ (see the display between (3.12) and (3.13) in \cite{V}).
\end{rmk}

\begin{rmk}\label{remtec}
In the proof of Theorem \ref{prelimdecay}, we are not able to ensure $\mu_n\to 0$, so we are compelled to distinguish cases \eqref{vanishing} and \eqref{lowerbound}. This basically relies in the type of estimate we want to prove (i.e., \eqref{decay}), which differ from \eqref{univest2}, that is, the natural one associated with a general nonlinearity. Indeed, \eqref{decay} has the same form of \eqref{univest1}, which works for power-type nonlinearities, whose scale invariance allows to get rid of the behavior of $\mu_n$. In detail, in \cite[page 13]{Souplet}, an argument by contradiction leads to
$$M_n(y_n)\geq M_n(x_n)>2n(1+d^{-1}(x_n,\partial\Omega_n))\geq 2n\to\infty\quad \mbox{as} \;\; n\to\infty,$$
which implies $\mu_n\to 0$, being $M_n$ a suitable positive power of $u_n$. Unlikely, in our setting, \eqref{a4} and \eqref{control} do not imply $\mu_n\to 0$.
\end{rmk}

Now we are ready to prove the main result of the paper.

\begin{proof}[Proof of Theorem \ref{mainthm}]
Take any $u\in\D^{1,p}_0(\R^N)$ solution to \eqref{prob}. According to Theorem \ref{regthm}, $u\in L^\infty(\R^N)\cap C^{1,\tau}_\loc(\R^N)$, since \eqref{hyp} guarantees $f\in L^\infty(\R^N)_+$.

Let us prove estimates \eqref{main1}--\eqref{main2}. Set $\kappa:=\frac{\kappa_0}{2}$, where $\kappa_0$ stems from Lemma \ref{univbound}. Then there exist $r',\kappa'>0$ (depending on $p,N,\Lambda,u$) such that $\|\nabla u\|_{L^p(B_{r'}^e)}\leq\kappa$ and $\|\nabla u\|_p\leq \kappa'$.

For any $R>0$ and $y\in\R^N$ we define
\begin{equation}
\label{blowup2}
u_R(y) := R^{\frac{N-p}{p-1}}u(Ry).
\end{equation}
A rescaling argument, together with \eqref{hyp}, yields
\begin{equation}
\label{scaledeq}
\begin{aligned}
-\Delta_pu_R &\leq R^N(f(Ry)+\Lambda u(Ry)^{p^*-1})\leq LR^N(1+R|y|)^{-\alpha}+\Lambda R^{-p'}u_R(y)^{p^*-1}\\
&\le LR^{N-\alpha}|y|^{-\alpha}+\Lambda R^{-p'}u_R(y)^{p^*-1} \quad \mbox{in} \;\; \R^N\setminus\{0\}.
\end{aligned}
\end{equation}
Take any $r>r'$ and consider $R\geq r$. Thus, writing $u_R^{p^*-1}=u_R^{p^*-p}u_R^{p-1}$ and applying Theorem \ref{prelimdecay} to $u$ entail
\begin{equation*}
\begin{aligned}
R^{-p'}u_R^{p^*-1}&\leq R^{-p'}\left(K_0 R^{\frac{N-p}{p-1}}(R|y|)^{\frac{p-N}{p}}\right)^{p^*-p}u_R^{p-1}\\
&= K_0^{p^*-p}|y|^{-p}u_R^{p-1}\leq K_0^{p^*-p} u_R^{p-1} \quad \mbox{for all} \;\; y\in B_1^e,
\end{aligned}
\end{equation*}
since $B_1^e\subseteq B_{r/R}^e$. Notice that $R^{N-\alpha}|y|^{-\alpha} \leq r^{N-\alpha}$ for all $y\in B_1^e$, owing to $\alpha>N$. Then we have
\begin{equation}
\label{locboundineq}
-\Delta_p u_R \leq a u_R^{p-1} + b \quad \mbox{in} \;\; B_1^e,
\end{equation}
being $a:=\Lambda K_0^{p^*-p}$ and $b:=Lr^{N-\alpha}$. Choose any $\eps>0$ such that $p+\eps<p_*$. Then, applying Theorem \ref{locboundthm} with $s:=p-1+\eps$, $\Omega:= B_9\setminus \overline{B}_1$, and $x\in\partial B_5$, we deduce
\begin{equation}
\label{supest}
\|u_R\|_{L^\infty(B_7\setminus B_3)} \leq C_\eps \left(\|u_R\|_{L^{p-1+\eps}(B_9\setminus B_1)}+1\right).
\end{equation}
for some $C_\eps>0$ depending on $p,N,\eps,r,\Lambda,L,\alpha$. Since $p-1+\eps<p_*-1$, \eqref{embedding} ensures
\begin{equation}
\label{weaklebbound}
\|u_R\|_{L^{p-1+\eps}(B_9\setminus B_1)} \leq C_\eps\|u_R\|_{L^{p_*-1,\infty}(B_9\setminus B_1)},
\end{equation}
enlarging $C_\eps$ if necessary. Hence, using \eqref{supest}--\eqref{weaklebbound} and Lemma \ref{weaklebboundlemma}, besides observing that $f\in L^1(\R^N)_+$ and $\|\cdot\|_{p_*-1,\infty}$ is invariant under the scaling \eqref{blowup2}, yield
\begin{equation}
\label{supestfinal}
\|u_R\|_{L^\infty(B_7\setminus B_3)} \leq \hat{C} \quad \mbox{for all} \;\; R\geq r.
\end{equation}
for a suitable $\hat{C}>0$ depending on $p,N,\Lambda,\alpha,u$. Therefore, \eqref{locboundineq} becomes 
$$-\Delta_p u_R \leq a C_1^{p-1} + b \quad \mbox{in} \;\; B_7\setminus B_3 $$
so that, enlarging $\hat{C}$ if necessary,
\begin{equation}
\label{supgradestfinal}
\|\nabla u_R\|_{L^\infty(B_6\setminus B_4)} \leq \hat{C} \quad \mbox{for all} \;\; R\geq r,
\end{equation}
owing to Lemma \ref{reglemma}. Finally, for any $x\in B_{5r}^e$, applying \eqref{supestfinal}--\eqref{supgradestfinal} with $R=\frac{|x|}{5}$ we obtain
\begin{equation*}
u(x)\leq C |x|^{\frac{p-N}{p-1}}\quad\mbox{and}\quad |\nabla u(x)|\leq C |x|^{\frac{1-N}{p-1}} \quad \mbox{for all} \;\; x\in B_{5r}^e,
\end{equation*}
where $C>0$ depends on $p,N,\Lambda,L,\alpha,r,u$. Exploiting Theorem \ref{regthm}, it turns out that
\begin{equation}
\label{finalests}
u(x)\leq C_1\left(1+|x|^{\frac{N-p}{p-1}}\right)^{-1} \quad \mbox{and} \quad |\nabla u(x)|\leq C_2\left(1+|x|^{\frac{N-1}{p-1}}\right)^{-1}
\end{equation}
for all $x\in\R^N$, for suitable $C_1,C_2>0$ depending on $p,N,\Lambda,L,\alpha,u$.

Now fix any $r>0$ and set $c:=\inf_{B_r} u$, which is positive by Theorem \ref{regthm}. Let $\Phi$ be the fundamental solution of the $p$-Laplacian (i.e., $\Delta_p \Phi=0$ in $\R^N\setminus\{0\}$) with $\Phi(|x|)=c$ for all $x\in \partial B_r$, that is,
\begin{equation}
\label{fundsol}
\Phi(x):=c\left(\frac{|x|}{r}\right)^{\frac{p-N}{p-1}} \quad \mbox{for all} \;\; x\in\R^N\setminus\{0\}.
\end{equation}
Then \cite[Proposition A.12]{GG} ensures $u\geq \Phi$ in $B_r^e$, so there exists $\hat{c}>0$ (depending on $p,N,r,c$) such that
$$ u(x) \geq \hat{c}|x|^{\frac{p-N}{p-1}} \quad \mbox{for all} \;\; x\in B_r^e. $$
Hence
\begin{equation}
\label{finalestu}
u(x)\geq C_0\left(1+|x|^{\frac{N-p}{p-1}}\right)^{-1} \quad \mbox{for all} \;\; x\in\R^N,
\end{equation}
where $C_0>0$ is an opportune constant depending on $p,N,u$.

Putting \eqref{finalests}--\eqref{finalestu} together entail \eqref{main1}--\eqref{main2}. It remains to prove \eqref{main3}. We proceed inspired by \cite{S16, EMSV}. Suppose by contradiction that \eqref{main3} is false: then there exists $\{x_n\}\subseteq \R^N$ such that $|x_n|\to+\infty$ and
\begin{equation}
\label{contradicthyp}
|x_n|^{\frac{N-1}{p-1}} |\nabla u(x_n)| \to 0 \quad \mbox{as} \;\; n\to\infty.
\end{equation}
For every $n\in\N$, set $R_n:=|x_n|$ and $u_n:=u_{R_n}$, where $u_{R_n}$ is defined through \eqref{blowup2}. Owing to \eqref{finalests}, $\{u_n\}$ is bounded in $L^\infty(\R^N)$, so that
\begin{equation}
\label{scaledeq2}
0\leq -\Delta_p u_n \leq LR_n^{N-\alpha}|x|^{-\alpha}+\Lambda R_n^{-p'}u_n^{p^*-1} \to 0 \quad \mbox{in} \;\; C^0_\loc(\R^N\setminus\{0\}),
\end{equation}
recalling \eqref{scaledeq}, $\alpha>N$, and $R_n\to+\infty$. Accordingly, Lemma \ref{reglemma} ensures that $\{u_n\}$ is bounded in $C^{1,\tau}_\loc(\R^N\setminus\{0\})$, so that $u_n\to u_\infty$ in $C^1_\loc(\R^N\setminus\{0\})$ for some $u_\infty\in C^1_\loc(\R^N\setminus\{0\})$, due to Ascoli-Arzelà's theorem. Moreover,
$$ \|\nabla u_n\|_p = \|\nabla u\|_p \quad \mbox{for all} \;\; n\in\N, $$
so that $u_n\rightharpoonup u_\infty$ in $\D^{1,p}_0(\R^N)$. Exploiting \eqref{scaledeq2}, we get
$$ -\Delta_p u_\infty = 0 \quad \mbox{in} \;\; \R^N\setminus\{0\}. $$
Then, either $u_\infty=0$ or $u_\infty=k\Phi$ for some $k>0$, where $\Phi$ is as in \eqref{fundsol}. According to \eqref{finalestu}, we have
$$ u_\infty(x) = \lim_{n\to\infty} u_n(x) \geq C_0 \left(1+|x|^{\frac{N-p}{p-1}}\right)^{-1} \quad \mbox{for all} \;\; x\in\R^N\setminus\{0\}, $$
so $u_\infty>0$ in $\R^N\setminus\{0\}$. Thus, $u_\infty=k\Phi$ for a suitable $k>0$. Set $y_n:=\frac{x_n}{R_n}\in \partial B_1$ for every $n\in\N$. Up to sub-sequences, $y_n\to \overline{y}$ for some $\overline{y}\in \partial B_1$. Via \eqref{contradicthyp} we obtain
$$ |\nabla u_n(y_n)| = |x_n|^{\frac{N-1}{p-1}}|\nabla u(x_n)| \to 0 \quad \mbox{as} \;\; n\to\infty. $$
Consequently,
\begin{equation*}
\begin{aligned}
|\nabla u_\infty(\overline{y})| &\leq |\nabla u_\infty(\overline{y})-\nabla u_\infty(y_n)|+|\nabla u_\infty(y_n)-\nabla u_n(y_n)|+|\nabla u_n(y_n)| \\
&\leq |\nabla u_\infty(\overline{y})-\nabla u_\infty(y_n)|+\|\nabla u_\infty-\nabla u_n\|_{C^0(\partial B_1)}+|\nabla u_n(y_n)| \to 0,
\end{aligned}
\end{equation*}
exploiting also $\nabla u_\infty\in C^0(\partial B_1)$ and $\nabla u_n\to \nabla u_\infty$ in $C^0(\partial B_1)$ as $n\to\infty$. This leads to $\nabla u_\infty(\overline{y})=0$, contradicting the fact that $u_\infty=k\Phi$ does not possess critical points (see \eqref{fundsol}). Hence \eqref{main3} is ensured, concluding the proof of Theorem \ref{mainthm}.
\end{proof}

\begin{rmk}
\label{unifests}
A careful inspection of the proofs reveals that the constants $C_1,C_2$ appearing in Theorem \ref{mainthm} depend on $N,p,\Lambda,L,\alpha$, as well as on $u$, via the following quantities: $\|\nabla u\|_p$, $\|u\|_{L^{p_*-1,\infty}(\R^N)}$, $\|u\|_{W^{1,\infty}(B_{5r})}$, and $\|\nabla u\|_{L^p(B_{r'}^e)}$, being $r>r'>0$ opportune. In particular, if no concentration of compactness occurs, a uniform (in $u$) control on \eqref{vanishsuperlevel}, \eqref{equiuniformint}, \eqref{layercakelimit}, and \eqref{reabsorbTh} can be provided. Indeed, we argue as follows.
\begin{itemize}
\item Through Sobolev's inequality, one has
\begin{equation}
\label{uniflimit}
|\Omega_k| \leq k^{-p^*} \int_{\Omega_k} u^{p^*} \dx \leq k^{-p^*} S^{-p^*/p}\|\nabla u\|_p^{p^*} \to 0 \quad \mbox{as} \;\; k\to+\infty,
\end{equation}
and the limit is uniform in $\|\nabla u\|_p$.
\item The lack of concentraction of compactness implies the $L^{p^*}$-equi-uniform integrability of $u$ (cfr. \cite[Lemma 3.2]{BG}), that is, for all $\eps>0$ there exists $\delta>0$ (independent of $u$) such that, for any measurable $\Omega\subseteq \R^N$ with $|\Omega|<\delta$,
$$ \int_\Omega u^{p^*} \dx < \eps. $$
Fix $\eps>0$. Reasoning as in \eqref{uniflimit},
$$ |\{u>M/2\}| \to 0 \quad \mbox{as} \;\; M\to+\infty $$
uniformly in $\|\nabla u\|_p$. Then we can choose $M>0$ such that $|\{u>M/2\}|<\delta$, where $\delta$ is provided by the equi-uniform integrability. Thus,
$$ \int_{\{u>M/2\}} (u-M/2)^{p^*} \dx \leq \int_{\{u>M/2\}} u^{p^*} \dx < \eps, $$
showing that
$$ \lim_{M\to+\infty} \int_{\{u>M/2\}} (u-M/2)^{p^*} \dx = 0 $$
uniformly in $\|\nabla u\|_p$, which ensures \eqref{equiuniformint}.
\item Owing to the layer cake lemma, \eqref{layercakelimit} is equivalent to
$$ \lim_{h\to 0} \int_0^h s^{p^*-1} |\{u>s\}| \ds = 0, $$
which now we have to guarantee uniformly in $u$. Let us fix any $\eps>0$. Since no concentration occurs then, up to sub-sequences (cfr. \cite[Lemma 3.2]{BG}),
\begin{equation}
\label{benaoum}
\lim_{R\to+\infty} \int_{B_R^e} u^{p^*} \dx = 0 \quad \mbox{uniformly in} \;\; u,
\end{equation}
so there exists $R>0$ (independent of $u$) such that $\int_{B_R^e} u^{p^*} \dx < \frac{p^*\eps}{2}$. On the other hand, let us choose $\delta>0$ such that $\delta<\left(\frac{p^*\eps}{2|B_R|}\right)^{\frac{1}{p^*}}$. Therefore, for all $h\in(0,\delta)$,
\begin{equation*}
\begin{aligned}
\int_0^h s^{p^*-1} |\{u>s\}| \ds& = \int_0^h s^{p^*-1} |\{u>s\}\cap B_R| \ds + \int_0^h s^{p^*-1} |\{u>s\}\cap B_R^e| \ds \\
&\leq |B_R|\frac{\delta^{p^*}}{p^*} + \int_0^{+\infty} s^{p^*-1}|\{u\chi_{B_R^e}>s\}| \ds \\&\qquad= \frac{1}{p^*}\left(|B_R|\delta^{p^*} + \|u\chi_{B_R^e}\|_{p^*}^{p^*}\right) < \frac{\eps}{2}+\frac{\eps}{2}=\eps,
\end{aligned}
\end{equation*}
where we used the layer cake lemma on $u\chi_{B_R^e}$. Since $R$ is independent of $u$, so is $\delta$, proving our claim.
\item To ensure \eqref{reabsorbTh}, we have to guarantee
\begin{equation}
\label{uniflimitTh}
\lim_{h\to 0} \int_{\R^N} T_h(u)^{p^*} \dx = 0 \quad \mbox{uniformly in} \;\; u.
\end{equation}
To this aim, fix any $\eps>0$. Due to \eqref{benaoum}, there exists $R>0$ such that $\int_{B_R^e} u^{p^*} \dx<\frac{\eps}{2}$. Choose any $\delta>0$ satisfying $\delta^{p^*}|B_R|<\frac{\eps}{2}$. Therefore, for all $h\in(0,\delta)$,
$$ \int_{\R^N} T_h(u)^{p^*} \dx \leq \int_{B_R} h^{p^*} \dx + \int_{B_R^e} u^{p^*} \dx = \delta^{p^*}|B_R|+\int_{B_R^e} u^{p^*} \dx < \frac{\eps}{2}+\frac{\eps}{2}=\eps. $$
Since $R>0$ does not depend on $u$, also $\delta$ is independent of $u$, ensuring \eqref{uniflimitTh}.
\end{itemize}
\end{rmk}

\section*{Acknowledgments}
The authors warmly thank Fabio De Regibus for the fruitful discussion about gradient estimates and symmetry results via the moving plane method.

The authors are member of the {\em Gruppo Nazionale per l'Analisi Ma\-te\-ma\-ti\-ca, la Probabilit\`a e le loro Applicazioni} (GNAMPA) of the {\em Istituto Nazionale di Alta Matematica} (INdAM). They are partially supported by the INdAM-GNAMPA Project 2024 titled {\em Regolarità ed esistenza per operatori anisotropi} (E5324001950001). \\
Laura Baldelli is partially supported by the ``Maria de Maeztu'' Excellence Unit IMAG, reference CEX2020-001105-M, funded by MCIN/AEI/10.13039/501100011033/. \\
This study was partly funded by: Research project of MIUR (Italian Ministry of Education, University and Research) PRIN 2022 {\em Nonlinear differential problems with applications to real phenomena} (Grant Number: 2022ZXZTN2).

\begin{small}

\end{small}
\end{document}